\documentclass[12pt,a4paper]{amsart}
\usepackage{amsmath,amssymb,amsfonts,amscd,latexsym,amsthm,mathrsfs}
\usepackage{color}
\usepackage{relsize}
\usepackage{hyperref}
\usepackage{fullpage}
\usepackage{enumitem}
\usepackage[top=2.0cm,left=2.0cm,right=2.0cm,bottom=2.0cm]{geometry}

\newtheorem{theorem}{Theorem}[section]
\newtheorem{lemma}[theorem]{Lemma}
\newtheorem{corollary}[theorem]{Corollary}

\theoremstyle{remark}

\newcommand{\C}{\mathbb C}
\newcommand{\R}{\mathbb R}

\newcommand{\Z}{\mathbb Z}
\newcommand{\N}{\mathbb N}

\newcommand{\hypt}
{{ {}_AF_B \left( {a_1,...,a_A \atop b_1,...,b_B }  \Bigm| t  \right) }}

\newcommand{\hyt}
{{ {}_0F_B \left( {* 
\atop b_1,...,b_B }  \Bigm| t  \right) }}

\title{An analogue of Siegel's determinant}
\author{By Tapani Matala-aho}
\address{Tapani Matala-aho, Aalto University, Department of Mathematics, P.O. Box 11100, FI-00076 Aalto, Finland }
\email{tapani.matala-aho@aalto.fi}

\begin{document}

\begin{abstract}
Siegel-Shidlovskii theory of $E$-functions involves a
non-vanishing proof for the determinants attached to
the linear forms $D^kR(t)$, derivatives
of an auxiliary function $R(t)$. 
Let a non-zero function $F(t)$ satisfy $m$th order
linear differential equation which 
we shall write using the differential operator $\Delta=tD$
and let $L(t)$ be any non-zero linear form
of the derivatives $\Delta^i F(t)$ $(i=0,...,m-1; m\ge 2)$.
The determinants $\det\mathcal A_k$ attached to the linear forms 
$\Delta^kL(t)$ have certain simple properties that allow
us to give a short proof for the non-vanishing of $\det\mathcal A_k$ for 
a class  of differential equations including a subclass of hypergeometric
differential equations.
\end{abstract}

\maketitle

\section{Introduction}

Let $\mathcal F$ be a differential field of characteristic zero,
say $\mathcal F=K((t))$,
where $K$ is the field of constants with respect
to the differential operator $D$.
If $K=\C$ or $K=\C_p$, then let $D=\frac{d}{dt}$ 
be the usual derivative operator.
In the following we suppose that a non-zero function 
$F:K\to K$ satisfies the linear homogeneous differential equation
\begin{equation}\label{10}
TD^mF(t)=Q_{1}D^{m-1}F(t)+...+Q_{m}D^0F(t)
\end{equation}
of order $m\ge2$, where 
$T, Q_{1},..., Q_{m}\in K[t]$ 
are polynomials with degrees
$S'=\deg T(t), r'_j=\deg Q_{j}(t)$.

In studying linear and algebraic independence of numbers,
Siegel \cite{Sie1}
 used Thue's lemma to construct an auxiliary function 
\begin{equation}\label{20}
R(t)=B_{0,1}D^{m-1}F(t)+...+B_{0,m}D^0F(t)
\end{equation}
with polynomials 
$B_{0,1},...,B_{0,m}\in K[t]$.
Siegel also introduced the following linear forms
\begin{equation}\label{30}
R_n=(TD)^n R(t)=
B_{n,1}D^{m-1}F(t)+...+B_{n,m}D^0F(t)
\quad \forall n\in\N,
\end{equation}
where 
$B_{n,1},...,B_{n,m}\in K[t]$.
A crucial tool in Siegel's method is a 
non-vanishing proof of the $m\times m$ determinants
\begin{equation}\label{40}
\det\mathcal B_k=
\left|\begin{array}{ccccc}
       & B_{k,1}     & ...  & B_{k,m}\\
       & .           &      & .  \\
       & .           &      & .  \\
       & .           &      & .  \\
       & B_{m+k-1,1} & ...  & B_{m+k-1,m}  \\
\end{array}
 \right|,\quad  k\in\N.
\end{equation}
Thereafter the approach shown by Siegel \cite{Sie1},\cite{Sie2} has widely been
applied  in proving algebraic independence of the values of
$E$-functions, see \cite{Shi}.
Shidlovskii \cite{Shi} made  major progress
by proving the non-vanishing of the determinants \eqref{40}
for a large class of functions including Siegel $E$-functions.
However, the non-vanishing proofs of the determinants \eqref{40} 
which use Shidlovskii's lemma have the disadvantage  
that they have not generally yielded  effective independence measures.
Further advance came from Beukers-Brownawell-Heckman \cite{BBH}
who proved Siegel's normality condition for a large class
of meromorphic hypergeometric functions 
\begin{equation}\label{50}
\hypt=\sum_{n=0}^{\infty}
 \frac{(a_1)_n\cdots(a_A)_n}{n!(b_1)_n\cdots (b_B)_n}t^n,
\end{equation}
thus giving a possibility to achieve effective quantitative
algebraic independence results for the values of the functions \eqref{50}.

Here we will use the differential operator
$\Delta=tD$ and write the differential equation \eqref{10} in an
equivalent form 
\begin{equation}\label{60}
N\Delta^mF(t)=P_{1}\Delta^{m-1}F(t)+...+P_{m}\Delta^0F(t)
\end{equation}
where $N, P_{1},..., P_{m}\in K[t]$ and $S=\deg N(t), r_j=\deg P_{j}(t)$.
Let
$$L(t)=
A_{0,1}\Delta ^{m-1}F(t)+...+A_{0,m}\Delta^0F(t),$$
where $A_{0,1},...,A_{0,m}\in K[t]$,
be an arbitrary linear form from which we will construct the following linear forms
\begin{equation}\label{70}
L_n=(N\Delta )^n L(t)=
A_{n,1}\Delta ^{m-1}F(t)+...+A_{n,m}\Delta^0F(t)
\quad\forall n\in\N,
\end{equation}
where again
$A_{n,1},...,A_{n,m}\in K[t]$.
So, let
\begin{equation}\label{80}
\det\mathcal A_k=
\left|
\begin{array}{ccccc}
       & A_{k,1}     & ...  & A_{k,m}\\
       & .           &      & .  \\
       & .           &      & .  \\
       & .           &      & .  \\
       & A_{k+m-1,1} & ...  & A_{k+m-1,m}  \\
\end{array}
\right|,\quad k\in\N,
\end{equation}
be the determinants analogous to Siegel's determinants \eqref{40}.
Under the conditions  
$$(A_{0,1},...,A_{0,m})\ne (0,...,0)$$
and 
$$r_m=S+1,\quad r_j< r_m\quad \forall j=1,...,m-1$$
we will prove in Theorem \ref{T1} the non-vanishing of the determinant 
$\det\mathcal A_k$ directly by studying the degrees of the
polynomials $A_{i,j}$. 
Consequently our approach will imply effective quantitative
linear independence results for a class of functions $F(t)$
satisfying \eqref{60}, \eqref{150} and \eqref{160} including e.g. the hypergeometric series
\begin{equation}\label{90}
\hyt=\sum_{n=0}^{\infty} \frac{t^n}{n!(b_1)_n\cdots (b_B)_n}.
\end{equation}
Because the determinant \eqref{80} is non-zero for any non-zero linear form $L(t)$,
we get immediately an irreducibility criterion, in Corollary \ref{C1},
for the differential operator \eqref{60}.

It is interesting to note that, if the method of Theorem \ref{T1}
is used with the operator $D$ and the corresponding
classical linear forms \eqref{30}, then we have the assumptions
of Theorem \ref{T2} which are similar to the assumption \eqref{160}. 
However, the assumptions used in Theorem \ref{T2}
do not correspond to those interesting cases, say \eqref{90}, which are covered by Theorem \ref{T1}. 

The method of Theorem \ref{T1} is applied also in \cite{AM} 
where some linear independence results are proved for the solutions
of certain $q$-functional equations analogous to \eqref{10}. 
Here, it is interesting that these solutions include $q$-analogues of the series \eqref{90}.

\section{The linear forms}

Let
\begin{equation}\label{100}
L_0=L(t)=
A_{0,1}\Delta ^{m-1}F(t)+...+A_{0,m}\Delta^0F(t)
\end{equation}
be any linear form, where $A_{0,1},...,A_{0,m}\in K[t]$.
Then we define the linear forms $L_n$ recursively by
\begin{equation}\label{110}
L_{n+1}=N\Delta L_n\quad\forall n\in\N
\end{equation}
and we denote
\begin{equation}\label{120}
L_n=A_{n,1}\Delta ^{m-1}F(t)+...+A_{n,m}\Delta^0F(t)
\quad\forall  n\in\N.
\end{equation}
The definitions \eqref{110} and \eqref{120} imply directly the following recurrences.
\begin{lemma}
The polynomials $A_{n,1}(t),..., A_{n,m}(t)$ satisfy the 
linear recurrences
\begin{equation}\label{130}
A_{k+1,j}=  P_{j}A_{k,1}+N\Delta A_{k,j}+NA_{k,j+1}
\quad\forall j=1,...,m-1, 
\end{equation}
\begin{equation}\label{140}
A_{k+1,m}=P_{m}A_{k,1}+N\Delta A_{k,m}
\end{equation}
for all $k\in\N$.
\end{lemma}

\section{ A non-vanishing proof }

Here we suppose that the linear form \eqref{100} is arbitrary.

\begin{theorem}\label{T1}
Let the linear form \eqref{100} satisfy
\begin{equation}\label{150}
(A_{0,1},...,A_{0,m})\ne (0,...,0)
\end{equation}
and assume that
\begin{equation}\label{160}
S\ge 0,\quad r_m=S+1,\quad r_j< r_m\quad \forall j=1,...,m-1.
\end{equation}
Then
\begin{equation}\label{170}
\det\mathcal A_k=\left|
\begin{array}{ccccc}
       & A_{k,1}     & ...  & A_{k,m}\\
       & .           &      & .  \\
       & .           &      & .  \\
       & .           &      & .  \\
       & A_{k+m-1,1} & ...  & A_{k+m-1,m}  \\
\end{array} 
\right|\ne 0,\quad \forall k\in\N.
\end{equation}
\end{theorem}

Proof. 
Denote $\deg A_{k,j}=d_j$,
and let us suppose
\begin{equation}\label{180}
d_1=A,...,d_{l-1}=B < d_{l}=C\ge D =d_{l+1},..., d_m=F 
\end{equation}
for some $1\le l\le m$ (if $l=1$, then $A=C$, and if $l=m$, then $F=C$).
(Here we put $\deg 0(x)=-\infty$, where $0(x)$ is the zero polynomial.)

Let first $l\ge 2$. Then
\begin{equation}\label{190}
\deg A_{k+1,l-1}=
\deg(P_{l-1} A_{k,1} +N \Delta A_{k,l-1}+NA_{k,l})=
\end{equation}
$$\max\{ r_{l-1}+d_1 ;\ S+d_{l-1} ;\ S+d_l \}= d_l + S.$$
If $h\le l-2$, then
\begin{equation}\label{200}
\deg A_{k+1,h}\le \max\{ r_{h}+d_1 ;\ S+d_{h} ;\ S+d_{h+1} \}\le d_l+S-1.
\end{equation}
If $l\le h\le m$, then
\begin{equation}\label{210}
\deg A_{k+1,h}\le 
\max\{ r_{h}+d_1 ;\ S+d_{h};\ S+d_{h+1}\}\le d_l+S.
\end{equation}
Let $l=1$, then
\begin{equation}\label{220}
\deg A_{k+1,h}\le
\max\{ r_{h}+d_1 ;\ S+d_{h} ;\ S+d_{h+1} \}\le d_1 + S < 
\end{equation}
$$d_1+S+1=\max\{r_{m}+d_1;\ S+d_{m}\}=\deg(P_{m} A_{k,1} +N \Delta A_{k,m})=
\deg A_{k+1,m},$$
for all $h\le m-1$.

Proceeding by induction we shall meet the situation
\begin{equation}\label{230}
\deg A_{k+l-1,1}\ge \deg A_{k+l-1,2},...,\deg A_{k+l-1,m}
\end{equation}
in the $l$th row of the determinant
\begin{equation}\label{240}
\det\mathcal A_k=\left|
\begin{array}{ccccccc}
& A_{k,1}\quad ...&A_{k,l-1}   & A_{k,l}  &A_{k,l+1} ... & A_{k,m}\\
& .            &A_{k+1,l-1}   & A_{k+1,l}  &          & .  \\
& .            &            &          &              & .  \\
& A_{k+l-1,1}  &            &          &              & .  \\
& .            &            &          &              & A_{k+l,m}  \\
& .            &            &          &              & .  \\
& A_{k+m-1,1}  &            & ...      &              & A_{k+m-1,m}  \\
\end{array}
 \right|.
\end{equation}

By \eqref{220} the process starts in a similar fashion 
from the $(l+1)$th row
and consequently the degrees will behave in the following manner
\begin{equation}\label{250}
\left|
\begin{array}{ccccccc}
& A       &...  &B       &< C   &\ge D\quad...&F\\
& .       &...  &  < C+S &\ge\ ...&           &A+S+1 \\
& .       &     &        &      &             &. \\
& C+(l-1)S&\ge  &...     &      &             &. \\
& A'      &     &        &      &             &< C+lS+1 \\
& .       &     &        &      &             &. \\
& .       &     &        &      &             &. \\
& .       &     &        &      &             &. \\
& .       &     &...     &      &< C+(m-1)S+1  \ge\ ...  &. \\
\end{array}
\right|.
\end{equation}

Next we shall permute the rows in the determinant 
$\det\mathcal A_k$ in such a way that we have the maximum degree
polynomials in the anti-diagonal
which shows that the degree of determinant $\det\mathcal A_k$ 
is 
\begin{equation}\label{260}
md_l+(m^2-m)S/2+m-l 
\end{equation}
and thus 
$$\det \mathcal A_k\ne 0.\qed$$

Next we use the differential equation satisfied by $F(t)$ in the form \eqref{10}.
Also here we suppose that the linear form \eqref{20} is arbitrary
i.e. without any assumption for the order of $R(t)$.
Then the method used in proof of Theorem \ref{T1} shows the following result.

\begin{theorem}\label{T2}
Let 
\begin{equation}\label{270}
(B_{0,1},...,B_{0,m})\ne (0,...,0),
\end{equation}
and
\begin{equation}\label{280}
S'\ge 0,\quad r'_m=S'+1,\quad r'_j< r'_m\quad \forall j=1,...,m-1
\end{equation}
or
\begin{equation}\label{290}
r'_m\ge 0,\quad S'=r'_m+1,\quad r'_j\le r'_m\quad \forall j=1,...,m. 
\end{equation}
Then
\begin{equation}\label{300}
\det \mathcal B_k\ne 0\quad\forall k\in\N.
\end{equation}
\end{theorem}

Note, that in the proof of Theorem \ref{T2} we use recurrences 
which we get just by replacing  the operator $\Delta$ by $D$ in 
\eqref{130} and \eqref{140}. Thus we may prove Theorem \ref{T2} by the assumption
\eqref{280} but a minor difference in the step which corresponds to \eqref{190} 
 gives the possibility for the extra assumption \eqref{290}.
  
\section{ Applications}

Let us recall some basic facts we need of homogeneous linear differential equations.
It is known that
\begin{equation}\label{350}
t^nD^n=\sum_{k=0}^{n}s_1(n,k)\Delta^k\quad\forall n\in\mathbb{N},
\end{equation}
\begin{equation}\label{351}
\Delta^n=\sum_{k=0}^{n}S_2(n,k)t^kD^k\quad\forall n\in\mathbb{N},
\end{equation}
where $s_1(n,k)$ and $S_2(n,k)$ are Stirling numbers of the first and second kind, respectively.  

Thus any linear homogeneous differential equation may be written as
\begin{equation}\label{310}
E_1(t,D)y=0
\end{equation}
or equivalently as
\begin{equation}\label{320}
E_2(t,\Delta)y=0,
\end{equation}
where $E_i(t,x)\in K[t,x]$ and $\deg_x E_1(t,x)=\deg_x E_2(t,x)$, 
the order of the differential equation.
Put $\Lambda_1=D,\Lambda_2=\Delta$ and consider the following conditions:

$1^*$ The differential equation 
\begin{equation}\label{321}
E_i(t,\Lambda_i)y=0
\end{equation}
for a particular solution $F$, is of the minimum order $m$, 
which is denoted by 
$$\deg_{\Lambda_i}(F)=m\quad i=1,2.$$

$2^*$ The derivatives
\begin{equation}\label{330}
\Lambda_i^0 F,...,\Lambda_i^{m-1}F
\end{equation}
of a particular solution $F$ of \eqref{321} are linearly independent over $K(t)$
for $i=1,2$.

$A^*$ The differential equation \eqref{321} is homogeneously linearly
irreducible i.e.
$$\deg_{\Lambda_i}(F)=m$$
for all solutions of \eqref{321} and $i=1,2$.

$B^*$ The derivatives
\begin{equation}\label{340}
\Lambda_i^0 F,...,\Lambda_i^{m-1}F
\end{equation}
are linearly independent over $K(t)$ for all solutions $F$ of \eqref{321} and $i=1,2$.

$C^*$ The differential operator $E(\Lambda_i,t)$ is irreducible for $i=1,2$.

The equivalence between the cases $i=1$ and $i=2$ follows from
the fact that the connection matrix $\mathcal S$ between the sets  \eqref{330}, $i=1,2$, 
 is subdiagonal having a non-zero determinant.
The cases $1^*$ and $2^*$ are equivalent and further
from Nesterenko (see \cite{Nes} Lemma 3.1.) it follows that that the conditions 
$A^*$, $B^*$ and $C^*$
are equivalent. Clearly $A^*$ implies $1^*$.

\begin{corollary}\label{C1}
Assume the condition \eqref{160}  of Theorem \ref{T1}.  Then the differential operator $E(\Delta ,t)$
corresponding to the equation \eqref{60} is irreducible.
\end{corollary}

Proof. Let $F$ be any non zero solution of \eqref{60}. Let us suppose on the contrary that
\begin{equation}\label{360}
L_0=A_{0,1}\Delta ^{m-1}F+...+A_{0,m}\Delta^0F=0
\end{equation}
with some $(A_{0,1},...,A_{0,m})\in K[t]^m\setminus\{(0,...,0)\}.$
Now \eqref{110} and \eqref{360} imply
\begin{equation}\label{370}
L_n=A_{n,1}\Delta ^{m-1}F+...+A_{n,m}\Delta^0F=0
\end{equation}
for all $n\in\N$ and thus
\begin{equation}\label{380}
\mathcal A_0 (\Delta^{m-1}F,...,\Delta^0F)^T=(0,...,0)^T.
\end{equation}
By \eqref{170} we have a
contradiction with the assumption that $F(t)$ is a non-zero function.

Thus  the differential operator $E(\Delta ,t)$ is irreducible
by the equivalence of   $B^*$ and $C^*$.\qed   

In a similar fashion we may prove a corresponding result for the differential operator $D$.

\begin{corollary}\label{Cor4.2}
Assume the condition \eqref{280} or \eqref{290} of Theorem \ref{T2}.  Then the differential operator $E(D ,t)$
corresponding to the equation \eqref{10} is irreducible.
\end{corollary}

On the other hand, due the folklore the above corollaries have converse statements.
Namely, if the differential operator is irreducible, then Corollary X implies Theorem X.
For the completeness a  sketch of the proof for this fact is given starting from the basic facts of differential modules.

Let 
$\mathcal{D}=\mathcal{F}[\partial]$
be a ring of differential operators, where, say, $\partial=D$ or $\partial=\Delta$.
If $L\in\mathcal D$, then $\mathcal DL$ denotes the left submodule of $\mathcal D$
generated by $L$. By these notations we may state
the following lemma which may be found e.g. from \cite{Sin}.
\begin{lemma} 
$L\in\mathcal D$ is irreducible
 $\Leftrightarrow$
$\mathcal D/\mathcal DL$ is a simple left $\mathcal D$-module.
\end{lemma}

Now we are ready to show  how Corollary \ref{C1} implies Theorem \ref{T1}.
Put 
$$W=\mathcal D/\mathcal DE(\Delta ,t)$$
and let $W^*$ be its dual.
Then 
$$m =dim_{K(t)} W=dim_{K(t)} W^*.$$ 
Define then $V^*=\mathcal DL$ be a left $\mathcal D$-submodule of $W^*$ 
generated by a non zero $L=L_0\in W^*$. 
Suppose now that $E(\Delta ,t)$ is irreducible. Hence
$W$ and consequently $W^*$ are simple implying that $V^*=W^*$.
But
$$\det \mathcal A_k=0 \Leftrightarrow \dim_{K(t)} V^*<m$$  
see, e.g. \cite{Shi}.
Thus the irreducibility of the differential operator $E(\Delta ,t)$
shows the non-vanishing of the determinant $\mathcal A_k$.\qed
 
\begin{theorem}
Theorem \ref{T1} is satisfied by the series 
\begin{equation}\label{390}
F(t)=\sum_{n=0}^{\infty} \frac{t^n}{\prod_{k=0}^{n-1} P(k)},\quad
\deg P(x)=m\ge 2,\quad P(-1)=0.
\end{equation}
\end{theorem}

Proof.
When one writes
$P(x)=(x+1)P_1(x+1)$, then the series \eqref{390}
satisfies the differential equation
\begin{equation}\label{400}
[\Delta P_1(\Delta)-t]F(t)=0
\end{equation}
where $S=r_j=0$ for all $j=1,...,m-1$ and $r_m=1$.
Thus Theorem \ref{T1} is valid for $F(t)$.\qed

Here we note that in Lemma 4.2 of \cite{BBH} the authors verified
the condition $2^*$ for the hypergeometric functions 
\begin{equation}\label{410}
\sum_{n=0}^{\infty} \frac{(a_1)_n\cdots(a_p)_n}
{n!(b_1)_n\cdots (b_{q-1})_n}t^n,\quad p<q,
\end{equation}
where $a_i-b_j\notin\Z$ for all $i,j$ with $b_q=1$.
Note also the similar results of Shalikhov, see e.g. \cite{Sal}.

\section{Appendix  by an anonymous referee}

\subsection{An irreducibility criterion.}
As in the introduction of the text, let $K = \mathbb C$ or $\mathbb C_p$, and let ${\mathcal K} = K((z))$ 
be the field of formal power series in the variable $z$, endowed with the derivations  
$$\partial = d/dz, \theta =z\partial.$$
Let further $v$ be the $z$-adic valuation on $\mathcal K$: for any non-zero $f \in \mathcal K$, $v(f)$ 
is the order of $f$ at $0$. Finally, let ${\mathcal D} = {\mathcal K}[\partial] = {\mathcal K}[\theta]$ be 
the ring of differential operators with coefficients in $\mathcal K$. A non zero element $L$ of $\mathcal D$ 
is called irreducible if it cannot be written as the product $L_1 L_2$ of two elements of $\mathcal D$ of orders strictly smaller than that of $L$.

Let 
$$L = q_m(z) \partial^m + ... + q_1(z) \partial + q_0(z) \in {\mathcal D}$$
 be a $m$-th order differential operator, with coefficients $q_m \neq 0, q_{m-1}, ..., q_0$ in ${\mathcal K}$. 
We can alternatively write $L$ in  terms of $\theta$, as follows;
 $$L = a_n(z) \theta^m + ... + a_1(z) \theta + a_0(z),$$
 where the coefficients  $a_j$ can be computed in terms of the coefficients $q_i,  i \geq j$. 
From the classical Formula \eqref{350} of the text and its ```converse" (expressing $t^nD^n$ in terms of the $\Delta^k$'s), one deduces that  
$$ \forall i = 0, ..., m : \min_{j = i, ..., m} v(a_j) = \min_{j = i, ..., m} (v(q_i) - i).$$
\medskip

The Newton polygon ${\mathcal N}(L)$ of $L$ is the convex hull of the set 
$\cup_{0 \leq i \leq m} \{(x,y) \in \R^2, 0 \leq x \leq i, y \geq v(q_i) -i\}$. 
By the preceding relations, it coincides with the  convex hull of the set 
$\cup_{0 \leq i \leq m} \{(x,y) \in \R^2, 0 \leq x \leq i, y \geq v(a_i)\}$. 
The slopes of its non-vertical  bordering lines are called the {\it slopes of $L$} (at the point $0$). 
In particular, for any non zero $q \in {\mathcal K}$, the differential operators $L$ and $qL$ have the same slopes. 
Clearly, the  slopes of $L$ are rational numbers, and we have:

\begin{lemma} (Katz \cite{Katz}).  Let $L \in {\mathcal D}$ be a differential operator of order $m > 0$. 
Assume that $L$ has only one slope $\lambda$, and that the exact denominator of the rational number  
$\lambda$ is  equal to $m$. Then, $L$ is irreducible  in the ring ${\mathcal D}$. 
\end{lemma}

{\em Proof}  Since $L$ is irreducible if and only if  ${\mathcal D}/{\mathcal D}L$ is a simple $\mathcal D$-module, 
this is just a rephrasing of the Irreducibility Criterion 2.2.8 of Katz \cite{Katz}. See also \cite{vdPS}, bottom of p.79.

\subsection{An alternative proof of Corollary \ref{C1}}

We here consider the field ${\bf K} =K(t)$, and recall the notations 
$$D = d/dt, \Delta = tD$$
 from the text. Setting $z = 1/t$, we can view ${\bf K}  $ as a subfield of ${\mathcal K} = K((z))$. 
The valuation induced  on $\bf K$ by the valuation $v$ of $\mathcal K$ then corresponds to the place $t = \infty$.  
In particular, for any polynomial $A \in K[t] \subset {\mathcal K}$, we have:
$$ v(A) = - \deg (A).$$
Consider the operator of order $m$
$$E(\Delta, t) = N(t)\Delta^m - P_{1}(t)\Delta^{m-1} - ...- P_{m}(t)\Delta^0$$
from Formula \eqref{60} of the text, and assume that its (polynomial) coefficients 
satisfy Condition \eqref{150} of Theorem \ref{T1}. Setting $z = 1/t$, so that $\theta = t d/dt = -z d/dz = -\Delta$, we have
 $E(\Delta, t) = L(\theta, z) = a_m(z) \theta^m + a_{m-1} \theta^{m-1} + ...  + a_0(z)$, where
 $$a_m(z) =   (-1)^mN(\frac{1}{z}), ~ a_{m-1}(z) = \pm P_1(\frac{1}{z}), ~..., ~ a_0(z) = \pm P_m(\frac{1}{z}),$$
In the notation of the text for the degrees of the polynomials $N, P_i$, we therefore have:
$$v(a_m) = -S,~ v(a_{m-1}) = -r_1, ~...,~ v(a_1) = -r_{m-1}, ~v(a_0) = -r_m.$$
Conditions \eqref{160} now read:
$$ v(a_0) = v(a_m) - 1, ~{\rm and}~\forall j = 1, ..., m-1,   v(a_{m-j}) > v(a_0) .$$
Therefore, all the  points $\big(x = j, y =  v(a_j)\big) \in \R^2,  j = 1, ..., m-1,$ lie on or above the line $y =  v(a_0) +1$, 
and the Newton polygon of $L = L(\theta,z)$ is bordered by a unique non-vertical line, joining the points $\big(0, v(a_0)\big)$ 
and $\big(m, v(a_m) = v(a_0) +1\big)$. In particular, ${\mathcal N}(L)$ has a unique slope 
$$\lambda = \frac{v(a_m ) - v(a_0)}{m} = \frac{1}{m}.$$
Since its denominator is  $m$, we deduce from the criterion that $L$ is irreducible in the ring ${\mathcal K}[\theta]$. 
A fortiori, $L$ is irreducible in the ring ${\bf K}[\theta] = K(t)[D]$, and so is $E = E(\Delta, t)$, as announced in Corollary \ref{C1}.

\subsection{An alternative proof of Corollary \ref{Cor4.2}}

With the same notations as above, let us now consider the operator of order $m$
$$E(D, t) = Q_0(t) D^m - Q_1(t)D^{m-1} -   ... - Q_m(t)D^0$$
from Formula \eqref{10} (where we have set : $Q_0(t) := T(t)$), and assume the degrees 
$r'_0 := S' = \deg(Q_0), r'_j = \deg(Q_j), j = 1, ...m$ of its (polynomial) coefficients satisfy either 
Condition \eqref{280} or Condition \eqref{290} of Theorem \ref{T2}. We shall then prove:

\begin{corollary} (= { \bf Corollary \ref{Cor4.2}}): under Condition \eqref{280} or \eqref{290}, the differential operator   
$E$ is irreducible in the ring ${\bf K}(t)[D]$.
\end{corollary}

Using Formula \eqref{350}, and then setting $z = 1/t$, we can rewrite
$$E(D,t) = \frac{Q_0(t)}{t^m}(t^mD^m) - \frac{Q_1(t)}{t^{m-1}}(t^{m-1}D^{m-1})  + ... - Q_m(t)$$
$$\quad  = \tilde Q_0 (t) \Delta^m + ~\tilde Q_1(t) \Delta^{m-1}+ ... + ~\tilde Q_m(t),$$
$$\qquad \qquad = \tilde a_m(z) \theta^m +  \tilde a_{m-1}(z) \theta^{m-1}+ ... + \tilde a_0(z) :=  \tilde L(\theta, z).$$
A computation similar to that of \S 1 (using both \eqref{350} and its converse) shows that the coefficients 
$\tilde a_j(z) = \pm \tilde Q_{m-j}(\frac{1}{z})$ of $\tilde L$ satisfy $v(\tilde a_m) = - \deg(Q_0) +m$, and more generally:
$$\forall i = 0, ..., m,   \min_{j=i,..., m} v(\tilde a_j) = \min_{j = i, ..., m} ( -\deg(Q_{m-j}) + j).$$
Therefore, the Newton polygon of $\tilde L$ is given by the convex hull of the points:
$$ \big(0,   -r'_m\big), \big(1, -r'_{m-1}+ 1\big),  ..., \big(m-1,  -r'_1 + m-1\big), \big(m,v(\tilde a_m) = -S' + m \big).$$

Let us now study this Newton polygon ${\mathcal N}(\tilde L)$.

{\it Case of Condition \eqref{280}} : in this case, the first point is $P_0 = \big(0,  -r'_m =  -(S'+1)\big)$, 
the last point is $P_m = \big(m,  -S' +m\big)$, and since $r'_i < r'_m$, i.e. $r'_i \leq S'$ for $1 \leq i \leq  m-1$, 
all the other points $\big(x= i, y = -r'_{m-i} + i \geq - S' + i \big)$ lie on or  above the line of slope 1 passing through 
the  point $P_m$, while $P_0$ lies  below this line. Therefore, ${\mathcal N}(\tilde L)$ has only one non vertical side, 
given by $[P_0P_m]$, whose slope is equal to 
$$ \lambda = \frac{(-S'+m) + (S'+1)}{m} = \frac{m+1}{m}.$$
 Since the exact denominator of this number is $m$,  the irreducibility Criterion applies.

{\it Case of Condition} \eqref{290} : in this case, the first point is $P_0 = \big(0, -r'_m = -(S'-1)\big)$, 
the last point is $P_m = \big(m,   -S' +m\big)$, the slope of the line $(P_0P_m)$ is
$$\lambda = \frac{(-S' +m) + (S'-1)}{m} = \frac{m-1}{m} = 1 - \frac{1}{m},$$
and since $r'_i \leq r'_m < S'$ for $1 \leq i \leq  m-1$, all the other points $P_i = \big( i,  -r'_{m-i} + i)$ 
lie above this line: indeed, the slope of $(P_iP_m)$ is given by
$$ \frac{(-S'+m) -(-r'_{m-i} +i)}{m-i} = 1 - \frac{S'-r'_{m-i}}{m-i} < 1 - \frac{1}{m}= \lambda.$$
Therefore, $\lambda$ is the unique slope of ${\mathcal N}(\tilde L)$, and since its  denominator is equal to $m$, 
the irreducibility Criterion again applies.

\subsection{Corollary X implies Theorem X} 

In each case, Corollary X says that a certain  $\mathcal D$-module $V$ is irreducible. 
In each case, Theorem X concerns a certain $\mathcal D$-module $W$, as well as the natural structure of 
$\mathcal D$-module that its dual $W^*$ acquires, and says the $\mathcal D$-submodule of $W^*$ generated by any  
non-zero linear form $\ell $ on $W$  fills up $W^*$. More precisely, the determinant considered 
in Formula \eqref{80} (resp. \eqref{40}) vanishes if and only if the dimension over $\mathcal K$ of the 
$\mathcal D$-submodule generated by $\ell = L$ (resp. $\ell = R$) is strictly smaller than $m =dim_{\mathcal K} W$.   
So, Theorem X exactly says that $W^*$ is irreducible (and this occurs if and only if $W$ itself is irreducible). 
To prove the equivalence of Corollary X and Theorem X, it thefore suffices to prove that $V$ is isomorphic to $W$ or to $W^*$. 
This follows from the fact that in each case $X = 1, resp.~ 2$, both $V$ and $W$ are isomorphic to 
${\mathcal D}/{\mathcal D}E$ (or to its dual), for the same differential operator $E$.

\end{document}